\documentclass{amsart}
\usepackage{amssymb}
\usepackage{upref}

\usepackage[all,cmtip]{xy}

\usepackage[lite,initials,shortalphabetic]{amsrefs}

\makeatletter
  \renewcommand{\p@enumi}{\thesubsection}
\makeatother

\newenvironment{resumeenumerate}[1]
{\begin{enumerate}
 \setcounter{enumi}{#1}
 \addtocounter{enumi}{-1}
}
{\end{enumerate}
}

\newenvironment{lettered}
{\begin{list}{\thelettercounter.}
 {\usecounter{lettercounter}\def\makelabel##1{\hss\llap{##1}}}
}
{\end{list}
}
\newcounter{lettercounter}
\renewcommand{\thelettercounter}{\Alph{lettercounter}}

\theoremstyle{plain}
\newtheorem*{thmA}{Theorem A}
\newtheorem*{thmB}{Theorem B}

\makeatletter
\newcommand{\emsection}[1]{%
  \par
  \addpenalty\@secpenalty
  \vskip 6 pt plus 9 pt
  \emph{#1.}\nobreak\enspace\ignorespaces
}
\makeatother

\newcommand{\intro}{%
  \goodbreak
  \vskip 6 pt plus 9 pt
}

\numberwithin{equation}{subsection}

\newcommand{\Comma}{\rlap{\enspace ,}}
\newcommand{\Period}{\rlap{\enspace .}}

\newcommand{\cat}[1]{\boldsymbol{#1}}

\newcommand{\RelCat}{\mathbf{RelCat}}
\newcommand{\RelnCat}{\mathbf{Rel}^{n}\mathbf{Cat}}
\newcommand{\RelnpoCat}{\mathbf{Rel}^{n+1}\mathbf{Cat}}
\newcommand{\Rk}{\mathbf{Rk}}
\newcommand{\Ry}{\mathbf{Ry}}

\newcommand{\Cat}{\mathbf{Cat}}

\newcommand{\simp}{\mathrm{s}}

\DeclareMathOperator{\colim}{colim}
\DeclareMathOperator{\map}{map}
\DeclareMathOperator{\Gr}{Gr}

\newcommand{\iso}{\approx}
\newcommand{\intersect}{\cap}

\newcommand{\rel}{_{\mathrm{rel}}}

\newcommand{\spacedcdots}{{\cdot\;\cdot\;\cdot}}

\newcommand{\adj}[4]{#1\negmedspace: #2\adjarrows #3:\negmedspace #4}
\newcommand{\adjarrows}{\mathchoice{\longleftrightarrow}
  {\leftrightarrow}
  {\leftrightarrow}
  {\leftrightarrow}}

\newcommand{\abs}[1]{\left |#1\right|}

\begin{document}

\title[$n$-relative categories]
{$n$-relative categories: {A} model for the homotopy theory of $n$-fold
  homotopy theories}

\author{C. Barwick}
\address{Department of Mathematics, Massachusetts Institute of
  Technology, Cambridge, MA 02139}
\email{clarkbar@math.mit.edu}

\author{D.~M. Kan}
\address{Department of Mathematics, Massachusetts Institute of
  Technology, Cambridge, MA 02139}

\date{December 14, 2010}
%\date{\today}

\begin{abstract}
  We introduce, for every integer $n \ge 1$, the notion of an
  \emph{$n$-relative category} and show that the category of the small
  $n$-relative categories is a model for the homotopy theory of
  \emph{$n$-fold homotopy theories}, i.e.\ homotopy theories of
  \ldots\ of homotopy theories.
\end{abstract}

\maketitle

%--------------------------------------------------------------------
%--------------------------------------------------------------------
\section{Background and motivation}
\label{sec:Bkgrnd}

In this introduction we
\begin{itemize}
\item recall some results of (higher) homotopy theory, and
\item explain how they led to the current manuscript.
\end{itemize}

\intro
We start with
%--------------------------------------------------------------------
\subsection{Rezk and re-Rezk}
\label{sec:reRezk}

In \cite{R} Charles Rezk constructed a left Bousfield localization of
the Reedy structure on the category $\simp\cat S$ of small simplicial
spaces (i.e.\ bisimplicial sets) and showed it to be a model for the
\emph{homotopy theory of homotopy theories}.

Furthermore it was noted in \cite{B} (and a proof thereof can be found
in \cite{Lu}*{\S1}) that iteration of Rezk's construction yields, for
every integer $n > 1$, a left Bousfield localization of the Reedy
structure on the category $\simp^{n}\cat S$ of small $n$-simplicial
spaces (i.e.\ $(n+1)$-simplicial sets) which is a model for the
\emph{homotopy theory of $n$-fold homotopy theories,} i.e.\ homotopy
theories of \ldots\ of homotopy theories.

We will call the weak equivalences in these left Bousfield
localization (which are often referred to as complete Segal
equivalences) just \emph{Rezk equivalences}.

\intro
Rezk's original result also gave rise to the following result on
%--------------------------------------------------------------------
\subsection{Relative categories}
\label{sec:RelCat}

Recall that a \emph{relative category} is a pair $(\cat C, w\cat C)$
consisting of a category $\cat C$ and a subcategory $w\cat C \subset
\cat C$ which contains all the objects of $\cat C$ and of which the
maps are called \emph{weak equivalences}.

Then it was shown in \cite{BK} that Rezk's model structure on
$\simp\cat S$ \eqref{sec:reRezk} can be lifted to a Quillen equivalent
\emph{Rezk structure} on the category $\RelCat$ of the small relative
categories, the weak equivalences of which will also
\eqref{sec:reRezk} be called \emph{Rezk equivalences}.

The category $\RelCat$ is connected to $\simp\cat S$ by a
\emph{simplicial nerve functor} $N\colon \RelCat \to \simp\cat S$ with
the property that a map $f \in \RelCat$ is a Rezk equivalence iff the
map $Nf \in\simp\cat S$ is so.  Moreover if we denote by $\Rk$ the
subcategories of the Rezk equivalences in both $\RelCat$ and
$\simp\cat S$, then the simplicial nerve functor has the property that
\begin{enumerate}
\item \label{relfnche}\emph{the relative functor}
  \begin{displaymath}
    N\colon (\RelCat, \Rk) \longrightarrow (\simp\cat S, \Rk)
  \end{displaymath}
  \emph{is a homotopy equivalence of relative categories}, in the
  sense that there exists a relative functor
  \begin{displaymath}
    M\colon (\simp\cat S, \Rk) \longrightarrow
    (\RelCat, \Rk)
  \end{displaymath}
  called a \emph{homotopy inverse} of $N$ such that the compositions
  $MN$ and $NM$ can be connected to the identity functors of $\RelCat$
  and $\simp\cat S$ by finite zigzags of natural weak equivalences.
\end{enumerate}

\intro
This in turn implies that
\begin{resumeenumerate}{2}
  \item \emph{the relative category $(\RelCat, \Rk)$ is, just like
      $(\simp\cat S, \Rk)$, a model for the homotopy theory of
      homotopy theories.}  
\end{resumeenumerate}

\intro
The proof of all this is essentially a relative version of the proof
of the following classical result of Bob Thomason.

%--------------------------------------------------------------------
\subsection{Thomason's result}
\label{sec:ThmRslt}

In \cite{T} Bob Thomason lifted the usual model structure on the
category $\cat S$ of small spaces (i.e.\ simplicial sets) to a Quillen
equivalent one on the category $\Cat$ of small categories and noted
that these two categories were connected by the \emph{nerve functor}
$N\colon \Cat \to \cat S$ which has the property that a map $f \in
\Cat$ is a weak equivalence iff $Nf \in \cat S$ is so.  It follows
that, if $\cat W$ denotes the categories of weak equivalences in both
$\Cat$ and $\cat S$, then
\begin{enumerate}
\item\label{HErf} \emph{the relative functor $N\colon (\Cat, \cat W)
    \to (\cat S, \cat W)$ is a homotopy equivalence of relative
    categories \eqref{relfnche}}
\end{enumerate}
which in turn implies that
\begin{resumeenumerate}{2}
  \item \label{RCmodl}\emph{the relative category $(\Cat, \cat W)$ is, just like
      $(\cat S, \cat W)$ a model for the theory of homotopy types.}
\end{resumeenumerate}

His proof was however far from simple as it involved notions like
two-fold subdivision and so-called Dwyer maps.

\intro
We end with recalling
%--------------------------------------------------------------------
\subsection{A result of Dana Latch}
\label{sec:RsltDLtch}

In \cite{La} Dana Latch noted that, if one just wanted to prove
\ref{HErf} and \ref{RCmodl}, one could do this by an argument that was
much simpler than Thomason's and that, instead of the cumbersome
two-fold subdivisions and Dwyer maps, involved the rather natural
notion of the \emph{category of simplices of a simplicial set}.

\intro
Now we can finally discuss
%--------------------------------------------------------------------
\subsection{The current paper}
\label{sec:CurPpr}

The results mentioned in \ref{sec:reRezk} and \ref{sec:RelCat} above
suggest that, for every integer $n > 1$, there might exist some
generalization of the notion of a relative category such that the
category of such generalized relative categories admits a model
structure which is Quillen equivalent to the Rezk structure on
$\simp^{n}\cat S$.

As however we did not see how to attack this question we turned to a
much simpler one suggested by the result of Dana Latch that was
mentioned in \ref{sec:RsltDLtch} above, namely to prove \ref{relfnche}
directly by showing that
\begin{itemize}
\item the \emph{simplicial nerve functor}
  \begin{displaymath}
    N\colon (\RelCat, \Rk) \longrightarrow (\simp\cat S, \Rk)
  \end{displaymath}
  has an appropriately defined \emph{relative category of bisimplices
    functor}
  \begin{displaymath}
    \Delta\rel\colon (\simp\cat S, \Rk) \longrightarrow
    (\RelCat, \Rk)
  \end{displaymath}
  as a homotopy inverse.
\end{itemize}

It turned out that not only could we do this, but the relative
simplicity of our proof suggested that a similar proof might work for
appropriately generalized relative categories.  And indeed, after the
necessary trial and error and frustration, we discovered a notion of
what we will call \emph{$n$-relative categories} which fitted the
bill.

Hence the current manuscript.

%--------------------------------------------------------------------
%--------------------------------------------------------------------
\section{An overview}
\label{sec:Ovrvw}

%--------------------------------------------------------------------
\subsection{Summary}
\label{sec:OvSmry}

There are five more sections.

\begin{itemize}
\item In the first (\S \ref{sec:nRlCt}) we introduce
  \emph{$n$-relative categories}.
\item In the second (\S \ref{sec:NSimpNerv}) we investigate an
  \emph{adjunction}
  \begin{displaymath}
    \adj{K}{\simp^{n}\cat S}{\RelnCat}{N}
  \end{displaymath}
  between the category $\simp^{n}\cat S$ of small $n$-simplicial
  spaces and the category $\RelnCat$ of small $n$-relative categories,
  in which the right adjoint $N$ is the \emph{$n$-simplicial nerve
    functor}.
\item Next (in \S \ref{sec:MainRlt} and \ref{sec:Prf2Prp}) we
  formulate and prove our main result.
\item In an appendix (\S \ref{sec:App}) we mention two relations
  between the categories $\RelnCat$ and $\RelnpoCat$.
\end{itemize}

\intro
In more detail:
%--------------------------------------------------------------------
\subsection{$n$-Relative categories}
\label{sec:nRelCat}

Motivated by the fact that in an $n$-simplicial space (i.e.\ an
$(n+1)$-simplicial set), just like in a simplicial space, the ``space
direction'' plays a different role than ``the $n$ simplicial
directions'', we define (in \S \ref{sec:nRlCt}) an \emph{$n$-relative
  category} $\cat C$ as an $(n+2)$ tuple
\begin{displaymath}
  \cat C = (a\cat C, v_{1}\cat C, \ldots, v_{n}\cat C, w\cat C)
\end{displaymath}
consisting of a category $a\cat C$ and subcategories $v_{1}\cat C$,
\ldots, $v_{n}\cat C$ and $w\cat C \subset a\cat C$ subject to the
following conditions:
\begin{enumerate}
\item \label{nRlCtcmm} Each of the subcategories contains all the
  objects of $a\cat C$ and together with $a\cat C$ they form a
  commutative diagram with $2n$ arrows of the form
  
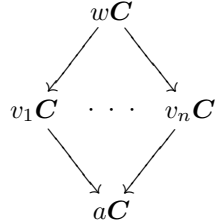
\begin{displaymath}
    \xymatrix@C=0.5em{
      & {w\cat C} \ar[dl] \ar[dr]\\
      {v_{1}\cat C} \ar[dr]
      & {\spacedcdots}
      & {v_{n}\cat C} \ar[dl]\\
      & {a\cat C}
    }
  \end{displaymath}
  which means that $v_{1}\cat C$, \ldots, $v_{n}\cat C$ can be
  considered as $n$ \emph{relative categories} which all have $w\cat
  C$ as their \emph{category of weak equivalences}.
\item \label{nRlCtamb} The purpose of $a\cat C$, the \emph{ambient
    category}, is to encode ``the extent to which any two of the
  $v_{i}\cat C$'s commute'' and we therefore impose on $a\cat C$ two
  conditions which simultaneously ensure that $a\cat C$ does not
  contain any superfluous information and that the associated (see
  \ref{sec:NSmpNrv}) $n$-simplicial nerve functor, just like the
  classical nerve functor, has a \emph{left adjoint} which is also a
  \emph{left inverse}.
\end{enumerate}

%--------------------------------------------------------------------
\subsection{The $n$-simplicial nerve functor}
\label{sec:NSmpNrv}

In \S \ref{sec:NSimpNerv} we introduce an adjunction
\begin{displaymath}
  \adj{K}{\simp^{n}\cat S}{\RelnCat}{N}
\end{displaymath}
between the category $\simp^{n}\cat S$ of the small $n$-simplicial
spaces (i.e.\ $(n+1)$-simplicial sets) and the category $\RelnCat$ of
the small $n$-relative categories \eqref{sec:nRelCat}.

If, for every integer $p \ge 0$, $\cat p$ denotes the category
\begin{displaymath}
  0 \longrightarrow \spacedcdots \longrightarrow p
\end{displaymath}
then the left adjoint $K$ is the colimit preserving functor which
sends each standard multisimplex $\Delta[p_{n}, \ldots, p_{1}, q]$ to
an $n$-relative version of the category
\begin{displaymath}
  \cat p_{n} \times \cdots \times \cat p_{1} \times \cat q
  \Period 
\end{displaymath}

The right adjoint $N$ will be referred to as the \emph{$n$-simplicial
  nerve functor}.

\intro
We also note that the counit and the unit of this adjunction have some
nice properties and in particular that \emph{the counit is an
  isomorphism} (which is equivalent to the statement that ``$K$ is not
only a left adjoint of $N$, but also a left inverse'' (cf.\
\ref{nRlCtamb})).

%--------------------------------------------------------------------
\subsection{The main result}
\label{sec:MnRslt}

To formulate our main result (in \S \ref{sec:MainRlt}) we use the
$n$-simplicial nerve functor $N$ \eqref{sec:NSmpNrv} to lift the Reedy
and the Rezk equivalences in $\simp^{n}\cat S$ \eqref{sec:reRezk} to
what we will also call Reedy and Rezk equivalences in $\RelnCat$ and
denote by $\Ry$ and $\Rk$ the subcategories of these Reedy and Rezk
equivalences in both $\simp^{n}\cat S$ and $\RelnCat$.

Our main result then is
\begin{thmA}
  The relative functor
  \begin{displaymath}
    N\colon (\RelnCat, \Rk) \longrightarrow (\simp^{n}\cat S, \Rk)
  \end{displaymath}
  is a homotopy equivalence of relative categories \eqref{relfnche}.
\end{thmA}

In view of the fact that the Rezk equivalences in $\simp^{n}\cat S$
are the weak equivalences in a left Bousfield localization of the
Reedy structure this theorem is a ready consequence of
\begin{thmB}
  The relative functor
  \begin{displaymath}
    N\colon (\RelnCat, \Ry) \longrightarrow (\simp^{n}\cat S, \Ry)
  \end{displaymath}
  is a homotopy equivalence of relative categories.
\end{thmB}

%--------------------------------------------------------------------
\subsection{The proof}
\label{sec:Prf}

Most of the proof of our main result \eqref{sec:MnRslt} is also in \S
\ref{sec:MainRlt}, except for the proof of two of the propositions
involved which we will deal with in \S \ref{sec:Prf2Prp}.

Apart from some properties of the counit and the unit of the
adjunction \eqref{sec:NSmpNrv}
\begin{displaymath}
  \adj{K}{\simp^{n}\cat S}{\RelnCat}{N} \Comma
\end{displaymath}
the proof involves the rather obvious \emph{category of multisimplices
functor}
\begin{displaymath}
  \Delta\colon \simp^{n}\cat S
    \longrightarrow \Cat \enspace\eqref{sec:ThmRslt}
\end{displaymath} 
and an $n$-relative version thereof, the \emph{$n$-relative category
  of multisimplices functor}
\begin{displaymath}
  \Delta\rel\colon \simp^{n}\cat S \longrightarrow
    \RelnCat \Comma
\end{displaymath}
which will be the required homotopy inverse \eqref{relfnche} of $N$.

In particular we need two rather simple properties of $\Delta$, as
well as two properties of $\Delta\rel$.  The proofs of the latter take
rather more effort and will therefore be dealt with separately in \S
\ref{sec:Prf2Prp}.

%--------------------------------------------------------------------
\subsection{An appendix}
\label{sec:Apdx}

In an appendix (\S \ref{sec:App}) we mention two relations between the
categories $\RelnCat$ and $\RelnpoCat$ which one would expect higher
homotopy theories to have:
\begin{em}
  \begin{enumerate}
  \item \label{ApxPrvi} That the functor $\RelnCat \to \RelnpoCat$
    which sends
    \begin{displaymath}
      (a\cat C, v_{1}\cat C, \ldots, v_{n}\cat C, w\cat C)
      \enspace\text{to}\enspace
      (a\cat C, v_{1}\cat C, \ldots, v_{n}\cat C, w\cat C, w\cat C)
    \end{displaymath}
    has a right adjoint which is a left inverse.
  \item \label{ApxPrvii} That every object of $\RelnpoCat$ gives rise
    to a category enriched over $\RelnCat$.
  \end{enumerate}
\end{em}

%--------------------------------------------------------------------
%--------------------------------------------------------------------
\section{$n$-Relative categories}
\label{sec:nRlCt}

After a brief review of \emph{relative categories} we
\begin{itemize}
\item introduce $n$-relative categories ($n \ge 1$) and
\item describe some simple but useful examples which we will need in
  the next section.
\end{itemize}

%--------------------------------------------------------------------
\subsection{Relative categories}
\label{sec:RlCt}

A \textbf{relative category} is a pair $(\cat C, \cat W)$ (often
denoted by just $\cat C$) consisting of a category $\cat C$ (the
\textbf{underlying category}) and a subcategory $\cat W \subset \cat
C$, the maps of which are called the \textbf{weak equivalences}, and
which is \emph{only} subject to the condition that it contains all the
\emph{objects} of $\cat C$ (and hence all the \emph{identity maps}).

The category of small relative categories and the \textbf{relative}
(i.e.\ weak equivalence preserving) functors between them will be
denoted by $\RelCat$.

Two relative functors $\cat C \to \cat D$ are called \textbf{naturally
  weakly equivalent} if they can be connected by a finite zigzag of
natural weak equivalences and a relative functor $f\colon \cat C \to
\cat D$ will be called a \textbf{homotopy equivalence} if there exists
a relative functor $g\colon \cat D \to \cat C$ (called a
\textbf{homotopy inverse} of $f$) such that the compositions $gf$ and
$fg$ are naturally weakly equivalent to $1_{\cat C}$ and $1_{\cat D}$
respectively.

%--------------------------------------------------------------------
\subsection{What to look for in a generalization}
\label{sec:LkGnrl}

In trying to generalize the notion of a relative category we were
looking for
\begin{itemize}
\item[$*$] a notion of $n$-relative category for which the associated
  $n$-simplicial nerve functor to $n$-simplicial spaces, just like the
  classical nerve functor, has a \emph{left adjoint} which is also a
  \emph{left inverse}.
\end{itemize}

Motivated by the fact that in an $n$-simplicial space (i.e.\ an
$(n+1)$-simplicial set), just like in a simplicial space, the ``space
direction'' plays a different role than ``the $n$ simplicial
directions'', we start with considering sequences
\begin{displaymath}
  \cat C = (a\cat C, v_{1}\cat C, \ldots, v_{n}\cat C, w\cat C)
  \enspace (n \ge 1)
\end{displaymath}
consisting of a category $a\cat C$ and subcategories $v_{1}\cat C$,
\ldots, $v_{n}\cat C$ and $w\cat C \subset a\cat C$, each of which
contains all the objects of $a\cat C$ and which together with $a\cat
C$ form a commutative diagram with $2n$ arrows of the form
\begin{displaymath}
  \xymatrix@C=0.5em{
    & {w\cat C} \ar[dl] \ar[dr]\\
    {v_{1}\cat C} \ar[dr]
    & {\spacedcdots}
    & {v_{n}\cat C} \ar[dl]\\
    & {a\cat C}
  }
\end{displaymath}

Such a sequence can be considered to consist of $n$ \emph{relative
  categories} $v_{1}\cat C$, \ldots, $v_{n}\cat C$ which each has the
same \emph{category of weak equivalences} $w\cat C$ and an
\emph{ambient category} $a\cat C$ which encodes the relations between
the $v_{i}\cat C$ ($1 \le i \le n$).

However the associated $n$-simplicial nerve functor
\eqref{sec:ReeRzEq} will only recognize those \emph{maps} in $a\cat C$
which are finite compositions of maps in the $v_{i}\cat C$ ($1 \le i
\le n$) and only those \emph{relations} which are a consequence of the
commutativity of those squares in $a\cat C$ which are of the form
\begin{displaymath}
  \xymatrix{
    {\cdot} \ar[r]^{x_{1}} \ar[d]_{y_{1}}
    & {\cdot} \ar[d]^{y_{2}}\\
    {\cdot} \ar[r]_{x_{2}}
    & {\cdot}
  }
\end{displaymath}
in which $x_{1}, x_{2} \in v_{i}\cat C$ and $y_{1}, y_{2} \in
v_{j}\cat C$ (where $i$ and $j$ are not necessarily distinct).

In order that the associated $n$-simplicial nerve functor has a left
inverse we therefore have to impose some restrictions on $a\cat C$ and
define as follows
%--------------------------------------------------------------------
\subsection{$n$-Relative categories}
\label{sec:NRlCt}

An \textbf{$n$-relative category} $\cat C$ will be an $(n+2)$-tuple
\begin{displaymath}
  \cat C = (a\cat C, v_{1}\cat C, \ldots, v_{n}\cat C, w\cat C)
\end{displaymath}
consisting of a category $a\cat C$ and subcategories
\begin{displaymath}
  v_{1}\cat C, \ldots, v_{n}\cat C
  \enspace\text{and}\enspace
  w\cat C \subset a\cat C
\end{displaymath}
each of which contains all the objects of $\cat C$ and which form a
commutative diagram with $2n$ arrows of the form
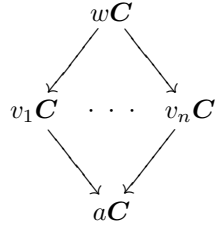
\begin{displaymath}
  \xymatrix@C=0.5em{
    & {w\cat C} \ar[dl] \ar[dr]\\
    {v_{1}\cat C} \ar[dr]
    & {\spacedcdots}
    & {v_{n}\cat C} \ar[dl]\\
    & {a\cat C}
  }
\end{displaymath}
and where $a\cat C$ is subject to the condition that
\begin{enumerate}
\item \label{nRlCtFcmp} every map in $a\cat C$ is a finite composition
  of maps in the $v_{i}\cat C$ ($1 \le i \le n$), and
\item \label{nRlCtCmt} every relation in $a\cat C$ is a consequence of
  the commutativity of those squares in $a\cat C$ which are of the
  form
  
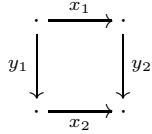
\begin{displaymath}
    \xymatrix{
      {\cdot} \ar[r]^{x_{1}} \ar[d]_{y_{1}}
      & {\cdot} \ar[d]^{y_{2}}\\
      {\cdot} \ar[r]_{x_{2}}
      & {\cdot}
    }
  \end{displaymath}
  in which $x_{1}, x_{2} \in v_{i}\cat C$ and $y_{1}, y_{2} \in
  v_{j}\cat C$ (where $i$ and $j$ are not necessarily distinct).
\end{enumerate}

%--------------------------------------------------------------------
\subsection{Some comments}
\label{sec:SmCmnt}

In an $n$-relative category $\cat C$, the categories $v_{1}\cat C$,
\ldots, $v_{n}\cat C$ are \emph{relative categories} which have $w\cat
C$ as their category of weak equivalences, and we will therefore
sometimes refer to the maps of $w\cat C$ as \textbf{weak
  equivalences}.

Moreover the category $a\cat C$ is more than a common underlying
category for the $v_{i}\cat C$ ($1 \le i \le n$) (as it may contain
additional relations) and will therefore be called the \textbf{ambient
category}.

Also note that
\begin{itemize}
\item [$*$] \emph{A $1$-relative category $\cat C$ is essentially just
  an ordinary relative category, as in that case $a\cat C = v_{1}\cat
  C$.}
\end{itemize}

%--------------------------------------------------------------------
\subsection{Relative functors}
\label{sec:RlFnc}

A \textbf{relative functor} $f\colon \cat C \to \cat D$ between two
$n$-relative categories $\cat C$ and $\cat D$ will be a functor
$f\colon a\cat C \to a\cat D$ such that
\begin{displaymath}
  fw\cat C \subset w\cat D
  \qquad\text{and}\qquad
  fv_{i}\cat C \subset v_{i}\cat D
  \qquad\text{for all $i \le i \le n$.}
\end{displaymath}

We will denote by $\RelnCat$ the resulting category of the small
$n$-relative categories and the relative functors between them.

\intro
We end with
%--------------------------------------------------------------------
\subsection{Some examples}
\label{sec:SmEx}

Some rather simple but useful examples of $n$-relative categories are
the following.

For every integer $p \ge 0$ let $\cat p$ denote the category
\begin{displaymath}
  0 \longrightarrow \spacedcdots \longrightarrow p
\end{displaymath}
and let $\abs{\cat p} \subset \cat p$ be its subcategory which
consists of the objects and their identity maps only.  Then we will
denote
\begin{enumerate}
\item \label{SmExi} by $\cat p^{w} \in \RelnCat$ the object such that
  \begin{displaymath}
    a\cat p^{w} = v_{i}\cat p^{w} = w\cat p^{w} = \cat p
    \qquad\text{for all $1 \le i \le n$}
  \end{displaymath}
\end{enumerate}
and
\begin{resumeenumerate}{2}
  \item \label{SmExii} by $\cat p^{v_{i}} \in \RelnCat$ the object
    such that
    \begin{displaymath}
      a\cat p^{v_{i}} = v_{i}\cat p^{v_{i}} = \cat p
      \qquad\text{and}\qquad
      v_{j}\cat p^{v_{i}} = w\cat p^{v_{i}} = \abs{\cat p}
      \qquad\text{for $j \ne i$}\Period
    \end{displaymath}
\end{resumeenumerate}

A simple calculation then yields that, for every sequence of integers
$p_{n}, \ldots, p_{1}, q \ge 0$,
\begin{resumeenumerate}{3}
  \item \label{SmExiii} $\cat p_{n}^{v_{n}}\times \cdots \times \cat
    p_{1}^{v_{1}} \in \RelnCat$ is such that
    \begin{align*}
      a(\cat p_{n}^{v_{n}}\times \cdots \times \cat p_{1}^{v_{1}})
      &= \cat p_{n}\times \cdots \times \cat p_{1}\\
      w(\cat p_{n}^{v_{n}}\times \cdots \times \cat p_{1}^{v_{1}})
      &= \abs{\cat p_{n}}\times \cdots \times \abs{\cat p_{1}}
      \qquad\text{and}\\
      v_{i}(\cat p_{n}^{v_{n}}\times \cdots \times \cat p_{1}^{v_{1}})
      &= \abs{\cat p_{n}}\times \cdots \times \cat p_{i}\times \cdots
      \times \abs{\cat p_{1}} \qquad\text{($1 \le i \le n$)}
    \end{align*}
\end{resumeenumerate}
and
\begin{resumeenumerate}{4}
  \item \label{SmExiv} $\cat p_{n}^{v_{n}}\times \cdots \times \cat
    p_{1}^{v_{1}} \times \cat q^{w} \in \RelnCat$ is such that
    \begin{align*}
      a(\cat p_{n}^{v_{n}}\times \cdots \times \cat p_{1}^{v_{1}}
      \times \cat q^{w})
      &= \cat p_{n}\times \cdots \times \cat p_{1}\times\cat q\\
      w(\cat p_{n}^{v_{n}}\times \cdots \times \cat p_{1}^{v_{1}}
      \times \cat q^{w})
      &= \abs{\cat p_{n}}\times \cdots \times \abs{\cat p_{1}}\times
      \cat q \qquad\text{and}\\
      v_{i}(\cat p_{n}^{v_{n}}\times \cdots \times \cat p_{1}^{v_{1}}
      \times\cat q^{w})
      &= \abs{\cat p_{n}}\times \cdots \times \cat p_{i}\times \cdots
      \times \abs{\cat p_{1}}\times \cat q \qquad\text{($1 \le i \le n$)}
    \end{align*}
\end{resumeenumerate}

%--------------------------------------------------------------------
%--------------------------------------------------------------------
\section{The $n$-simplicial nerve functor}
\label{sec:NSimpNerv}

We now
\begin{itemize}
\item introduce an adjunction (\ref{sec:reRezk} and \ref{sec:RlFnc})
  \begin{displaymath}
    \adj{K}{\simp^{n}\cat S}{\RelnCat}{N}
  \end{displaymath}
  in which the right adjoint $N$ is the \emph{$n$-simplicial nerve
    functor} which we mentioned in \ref{sec:LkGnrl},
\item use $N$ to lift the Reedy and the Rezk equivalences
  \eqref{sec:reRezk} from $\simp^{n}\cat S$ to $\RelnCat$,
\item note that the unit of the above adjunction has two nice
  properties, and
\item note that the counit is an isomorphism which is the same as
  saying that $N$ has a \emph{left adjoint} which is also a \emph{left
    inverse} (cf.\ \ref{sec:LkGnrl}).
\end{itemize}

\intro
We start with
%--------------------------------------------------------------------
\subsection{The adjunction $\adj{K}{\simp^{n}\cat S}{\RelnCat}{N}$}
\label{sec:Adjunction}

The \textbf{$n$-simplicial nerve functor} will be the right adjoint in
the adjunction (\ref{sec:reRezk} and \ref{sec:RlFnc})
\begin{displaymath}
  \adj{K}{\simp^{n}\cat S}{\RelnCat}{N}
\end{displaymath}
in which
\begin{enumerate}
\item \label{Adji} $N$ sends an object $\cat C \in \RelnCat$ to the
  $(n+1)$-simplicial set which as as its $(p_{n}, \ldots, p_{1},
  q)$-simplices $(p_{n}, \ldots, p_{1}, q \ge 0)$ the maps
  \eqref{sec:SmEx}
  \begin{displaymath}
    \cat p_{n}^{v_{n}}\times \cdots \times \cat p_{1}^{v_{1}} \times
    \cat q^{w} \longrightarrow \cat C \in \RelnCat
  \end{displaymath}
\end{enumerate}
and
\begin{resumeenumerate}{2}
  \item \label{Adjii} $K$ is the colimit preserving functor which, for
    every $n+1$ integers $p_{n}, \ldots, p_{1}, q \ge 0$, sends the
    \emph{standard $(p_{n}, \ldots, p_{1}, q)$-simplex} $\Delta[p_{n},
    \ldots, p_{1}, q]$ to
    \begin{displaymath}
      \cat p_{n}^{v_{n}}\times \cdots \times \cat p_{1}^{v_{1}} \times
      \cat q^{w} \in \RelnCat \Period
    \end{displaymath}
\end{resumeenumerate}

\intro
Using the functor $N$ we then can define
%--------------------------------------------------------------------
\subsection{Reedy and Rezk equivalences in $\RelnCat$}
\label{sec:ReeRzEq}

A map $f \in \RelnCat$ will be called a \textbf{Reedy} or a
\textbf{Rezk equivalence} if the map $Nf \in \simp^{n}\cat S$ is so
\eqref{sec:reRezk}, and we will denote by
\begin{displaymath}
  \text{$\Ry$ and $\Rk \subset \simp^{n}\cat S$}
  \qquad\text{and}\qquad
  \text{$\Ry$ and $\Rk \subset \RelnCat$}
\end{displaymath}
the subcategories of the Reedy and the Rezk equivalences in both
$\simp^{n}\cat S$ and $\RelnCat$.

\intro
Next we note a very useful property of the functors $K$ and $N$.
%--------------------------------------------------------------------
\subsection{The $2$-skeleton property}
\label{sec:2skel}
\begin{em}
  \begin{enumerate}
  \item \label{sklprpi} For every object $X \in \simp^{n}\cat S$, the
    $n$-relative category $KX$ is completely determined by the
    $2$-skeleton of $X$, i.e.\ the smallest subobject that contains
    all its multisimplices of total dimension $\le 2$, and
  \item \label{sklprpii} for every object $\cat C \in \RelnCat$, the
    $(n+1)$-simplicial set $N\cat C$ is completely determined by its
    $2$-skeleton and is in fact its own $2$-coskeleton.
  \end{enumerate}
\end{em}

\emsection{Proof}
This follows by a straightforward calculation from the observation
that the category \eqref{sec:SmEx}
\begin{displaymath}
  a(\cat p_{n}^{v_{n}}\times \cdots\times \cat p_{1}^{v_{1}}\times
  \cat q^{w}) = \cat p_{n}\times \cdots \times \cat p_{1}\times \cat q
\end{displaymath}
is a poset which has an \emph{object}, a \emph{generating map} and a
\emph{relation} for every multisimplex of the standard $(p_{n},
\ldots, p_{1}, q)$-simplex $\Delta[p_{n}, \ldots, p_{1}, q]$ in total
dimensions $0$, $1$ and $2$ respectively.

\intro
The $2$-skeleton property \eqref{sec:2skel} readily implies
%--------------------------------------------------------------------
\subsection{Proposition}
\label{sec:2skeliso}
\begin{em}
  For every object $\cat C \in \RelnCat$, the counit map
  \begin{displaymath}
    \varepsilon\cat C\colon KN\cat C \longrightarrow
    \cat C \in \RelnCat
  \end{displaymath}
  is an isomorphism.
\end{em}

\intro
Moreover, in view of the fact that the composition
\begin{displaymath}
  N\cat C \xrightarrow{\quad \eta N\cat C\quad}
  NKN\cat C \xrightarrow{\quad N\varepsilon\cat C\quad}
  N\cat C \in \simp^{n}\cat S
\end{displaymath}
is the identity, \ref{sec:2skeliso} implies
%--------------------------------------------------------------------
\subsection{Proposition}
\label{sec:PrpUnIso}
\begin{em}
  For every object $\cat C \in \RelnCat$, the unit map
  \begin{displaymath}
    \eta N\cat C\colon N\cat C \longrightarrow NKN\cat C \in
    \simp^{n}\cat S
  \end{displaymath}
  is an isomorphism.
\end{em}

\intro
We also note
%--------------------------------------------------------------------
\subsection{Proposition}
\label{sec:PrpUnReEq}
For every standard multisimplex $\Delta[p_{n}, \ldots, p_{1}, q] \in
\simp^{n}\cat S$, the unit map
\begin{displaymath}
  \eta\Delta[p_{n}, \ldots, p_{1}, q]\colon \Delta[p_{n}, \ldots,
  p_{1}, q] \longrightarrow NK\Delta[p_{n}, \ldots, p_{1}, q] \in
  \simp^{n}\cat S
\end{displaymath}
is a Reedy equivalence.

\emsection{Proof}
Note that
\begin{displaymath}
  \Delta[p_{n}, \ldots, p_{1}, q]
  = \Delta[p_{n},=]\times \cdots\times\Delta[=, p_{i},=]\times
  \cdots \times\Delta[=,q]
\end{displaymath}
where the $=$'s denote sequences of $0$'s and that
\begin{displaymath}
  K\Delta[p_{n},\ldots, p_{1},q] =
  \cat p_{n}^{v_{n}}\times \cdots \times \cat p_{n}^{v_{1}}\times \cat
  q^{w} \Period
\end{displaymath}
A straightforward calculation then yields that
\begin{displaymath}
  \eta\Delta[=,p_{i},=]\colon \Delta[=,p_{i},=] \longrightarrow
  NK\Delta[=,p_{i},=] = N\cat p_{i}^{v_{i}} \in \simp^{n}\cat S
\end{displaymath}
is an isomorphism for all $1 \le i \le n$, and that
\begin{displaymath}
  \eta\Delta[=,q]\colon \Delta[=,q] \longrightarrow
  NK\Delta[=,q] \in \simp^{n}\cat S
\end{displaymath}
is a Reedy equivalence, and the desired result now follows from the
fact that $N$ is a right adjoint and hence preserves products.

%--------------------------------------------------------------------
%--------------------------------------------------------------------
\section{The main result}
\label{sec:MainRlt}

Now we are ready for
\begin{itemize}
\item our main result, and
\item a proof thereof, except for the verification of two propositions
  which we put off till \S \ref{sec:Prf2Prp}.
\end{itemize}

\intro
We thus start with stating
%--------------------------------------------------------------------
\subsection{Theorem}
\label{sec:ThmRlFnHE}
\begin{enumerate}
\item \label{MnRslti} \emph{The relative functor}
  (\ref{sec:Adjunction} and \ref{sec:ReeRzEq})
  \begin{displaymath}
    N\colon (\RelnCat, \Rk)\longrightarrow (\simp^{n}\cat S, \Rk)
  \end{displaymath}
  \emph{is a homotopy equivalence \eqref{sec:RlCt}, and hence}
\item \label{MnRsltii} \emph{the relative category $(\RelnCat, Rk)$
    is, just like $(\simp^{n}\cat S, \Rk)$ \eqref{sec:reRezk}, a model
  for the homotopy theory of $n$-fold homotopy theories.}
\end{enumerate}

\intro
To prove this, it suffices, in view of the fact that the Rezk
equivalences in $\simp^{n}\cat S$ are the weak equivalences in a left
Bousfield localization of the Reedy structure, to show
%--------------------------------------------------------------------
\subsection{Theorem}
\label{sec:ThmRelFuncHE}
\begin{em}
  The relative functor \eqref{sec:ReeRzEq}
  \begin{displaymath}
    N\colon (\RelnCat, \Ry) \longrightarrow (\simp^{n}\cat S, \Ry)
  \end{displaymath}
  is a homotopy equivalence \eqref{sec:RlCt}.
\end{em}

\intro
In preparation for a proof we first discuss
%--------------------------------------------------------------------
\subsection{The category of multisimplices}
\label{sec:CatMltsmpl}

Let $\Delta[-] \subset \simp^{n}\cat S$ denote the full subcategory
spanned by the standard multisimplices.

Given an object $X \in \simp^{n}\cat S$, one then defines its
\textbf{category of multisimplices} $\Delta X$ as the over category
\begin{displaymath}
  \Delta X = \Delta[-]\downarrow X \Period
\end{displaymath}
Clearly this category is \emph{natural in $X$}.  Moreover it comes
with a \emph{forgetful functor}
\begin{displaymath}
  F\colon \Delta X \longrightarrow \simp^{n}\cat S
\end{displaymath}
which sends an object $\Delta[p_{n}, \ldots, p_{1}, q] \to X$ to the
object $\Delta[p_{n}, \ldots, p_{1}, q] \in \simp^{n}\cat S$.

One then readily verifies that, as in the classical case, the
resulting $\Delta X$-diagram of standard multisimplices has the
following properties.
%--------------------------------------------------------------------
\subsection{Proposition}
\label{sec:PrpClmSmp}
\begin{em}
  For every object $X \in \simp^{n}\cat S$, the obvious map
  \begin{displaymath}
    \colim_{\Delta X} F \longrightarrow X \in \simp^{n}\cat S
  \end{displaymath}
  is an isomorphism.
\end{em}

%--------------------------------------------------------------------
\subsection{Proposition}
\label{sec:PrpSmpRee}
\emph{For every object $X \in \simp^{n}\cat S$, the category $\Delta
  X$ is a Reedy category with fibrant constants
  \cite{H}*{15.10.1(2)}}.

\intro
Next we introduce an $n$-relative version $\Delta\rel$ of the above
functor $\Delta$, which is the prospective homotopy inverse of the
$n$-simplicial nerve functor.
%--------------------------------------------------------------------
\subsection{The $n$-relative category of multisimplices functor}
\label{sec:NRlFnc}

Let $\Delta\rel[-]$ denote the $n$-relative category such that
\begin{enumerate}
\item \label{NRli} $a\Delta\rel[-] = \Delta[-]$
  \eqref{sec:CatMltsmpl}.
\item \label{NRlii} $v_{i}\Delta\rel[-]$ ($1 \le i \le n$) is the
  subcategory of $\Delta[-]$ consisting of the maps
  \begin{displaymath}
    \Delta[p_{n}, \ldots, p_{1}, q] \longrightarrow
    \Delta[p_{n}', \ldots, p_{1}', q'] \in \simp^{n}\cat S
  \end{displaymath}
  for which the associated map $\cat p_{i} \to \cat p_{i}'$ sends the
  object $p_{i} \in \cat p_{i}$ to the object $p_{i}' \in \cat
  p_{i}'$, and
\item \label{NRliii} $w\Delta\rel[-] = v_{1}\Delta\rel[-] \intersect
  \cdots \intersect v_{n}\Delta\rel[-]$.
\end{enumerate}

Given an object $X \in \simp^{n}\cat S$ we then define its
\textbf{$n$-relative category of multisimplices} $\Delta\rel X$ as the
$n$-relative over category
\begin{displaymath}
  \Delta\rel X = \Delta\rel[-] \downarrow X \Period
\end{displaymath}

Clearly $\Delta\rel X$ is \emph{natural in $X$} and the resulting
functor
\begin{displaymath}
  \Delta\rel\colon \simp^{n}\cat S \longrightarrow \RelnCat
\end{displaymath}
has the following two properties which we will need in the proof of
\ref{sec:ThmRelFuncHE} but which will only be proved in
\ref{sec:Prf5p7}--\ref{sec:PrpAmbiso} and
\ref{sec:Prf5p8}--\ref{sec:PrpCanIso} below respectively.
%--------------------------------------------------------------------
\subsection{Proposition}
\label{sec:ObvIso}

\begin{em}
  For every object $X \in \simp^{n}\cat S$, the obvious maps
  \eqref{sec:CatMltsmpl}
  \begin{align*}
    \colim_{\Delta X}\Delta\rel F &\longrightarrow
    \Delta\rel X \in \RelnCat, \qquad\text{and}\\
    \colim_{\Delta X}N\Delta\rel F &\longrightarrow
    N\Delta\rel X \in \simp^{n}\cat S
  \end{align*}
  are isomorphisms.
\end{em}

%--------------------------------------------------------------------
\subsection{Proposition}
\label{sec:PrpNTrReEq}
\begin{em}
  There exists a natural transformation \eqref{sec:Adjunction}
  \begin{displaymath}
    \pi_{t}\colon \Delta\rel \longrightarrow K
  \end{displaymath}
  with the property that, for every standard multisimplex
  $\Delta[p_{n}, \ldots, p_{1}, q]$, the map
  \begin{displaymath}
    N\pi_{t}\Delta[p_{n}, \ldots, p_{1}, q]\colon
    N\Delta\rel\Delta[p_{n}, \ldots, p_{1}, q] \longrightarrow
    NK\Delta[p_{n}, \ldots, p_{1}, q] \in \simp^{n}\cat S
  \end{displaymath}
  is a Reedy equivalence.
\end{em}

\intro
Now we are ready for
%--------------------------------------------------------------------
\subsection{A proof of theorem \ref{sec:ThmRelFuncHE}}
\label{sec:PrfTh5p2}

To prove that the functors $N\Delta\rel$ and $1_{\simp^{n}\cat S}$ are
naturally Reedy equivalent we consider, for every object $X \in
\simp^{n}\cat S$, the commutative diagram (\ref{sec:CatMltsmpl},
\ref{sec:NRlFnc} and \ref{sec:PrpNTrReEq})
\begin{displaymath}
  \xymatrix@C=4em{
    {\colim_{\Delta X}N\Delta\rel F} \ar[r]^-{N\pi_{t}F} \ar[d]
    & {\colim_{\Delta X}NKF} \ar[d]
    & {\colim_{\Delta X}F} \ar[l]_-{\eta F}\ar[d]\\
    {N\Delta\rel X} \ar[r]^-{N\pi_{t}X}
    & {NKX}
    & {X} \ar[l]_-{\eta X}
  }
\end{displaymath}
in which the vertical maps are the obvious ones.

The vertical maps on the outside are, in view of \ref{sec:PrpClmSmp}
and \ref{sec:ObvIso}, isomorphisms and it thus suffices to prove that
the upper maps are Reedy equivalences.  But this follows immediately
from \ref{sec:PrpNTrReEq} and \ref{sec:PrpUnReEq}, and
\ref{sec:PrpSmpRee} and the result \cite{H}*{15.10.9(2)} that the
colimit of an objectwise weak equivalence between Reedy cofibrant
diagrams indexed by a Reedy category with fibrant constants is also a
weak equivalence.

\intro
Note that the fact that these four maps are Reedy equivalences also
implies that the functor $N\Delta\rel$ preserves Reedy equivalences
and so does therefore \eqref{sec:ReeRzEq} the functor $\Delta\rel$.

To prove that the functors $\Delta\rel N$ and $1_{\RelnCat}$ are also
naturally Reedy equivalent it suffices to show that, for every object
$\cat C \in \RelnCat$, both maps in the sequence
\begin{displaymath}
  \Delta\rel N\cat C \xrightarrow{\quad \pi_{t}N\cat C\quad}
  KN\cat C \xrightarrow{\quad \varepsilon\cat C\quad}
  \cat C \in \RelnCat
\end{displaymath}
are Reedy equivalences.  For the second map this follows from
\ref{sec:2skeliso}.  To deal with the first one we have to show that
$N\pi_{t}N\cat C$ is a Reedy equivalence in $\simp^{n}\cat S$.  This
we do by considering the above diagram for $X \in N\cat C$
\begin{displaymath}
  \xymatrix@C=4em{
    {\colim_{\Delta N\cat C}N\Delta\rel F} \ar[r]^-{N\pi_{t}F} \ar[d]
    & {\colim_{\Delta N\cat C}NKF} \ar[d]
    & {\colim_{\Delta N\cat C}F} \ar[l]_-{\eta F}\ar[d]\\
    {N\Delta\rel N\cat C} \ar[r]^-{N\pi_{t}N\cat C}
    & {NKN\cat C}
    & {N\cat C} \ar[l]_-{\eta N\cat C}
  }
\end{displaymath}
and then noting that all its maps are Reedy equivalences in view of
the fact that
\begin{enumerate}
\item \label{Fcti} the upper and the outside vertical maps are so by
  the above,
\item \label{Fctii} the map $\eta N\cat C$ is so in view of
  \ref{sec:PrpUnIso}, and
\item \label{Fctiii} Reedy equivalences have the two out of three
  property.
\end{enumerate}

%--------------------------------------------------------------------
%--------------------------------------------------------------------
\section{A proof of propositions \ref{sec:ObvIso} and \ref{sec:PrpNTrReEq}}
\label{sec:Prf2Prp}

We start with
%--------------------------------------------------------------------
\subsection{A proof of \ref{sec:ObvIso}}
\label{sec:Prf5p7}

Proposition~\ref{sec:ObvIso} is a ready consequence of the more
general proposition \ref{sec:PrpAmbiso} below.  To formulate the
latter we will, for every pair of objects $\cat C, \cat D \in
\RelnCat$, denote by $\map(\cat C, \cat D)$ the \emph{set} of maps
$\cat C \to \cat D \in \RelnCat$.  Then we can state
%--------------------------------------------------------------------
\subsection{Proposition}
\label{sec:PrpAmbiso}

\begin{em}
  Let $\cat T \in \RelnCat$ have an ambient category which is a poset
  with a terminal object $T$.  Then, for every object $X \in
  \simp^{n}\cat S$, the obvious map \eqref{sec:CatMltsmpl}
  \begin{displaymath}
    \colim_{\Delta X}\map(\cat T, \Delta\rel F) \longrightarrow
    \map(\cat T, \Delta\rel X)
  \end{displaymath}
  is an isomorphism.
\end{em}

\emsection{Proof}
One readily verifies that the given map is onto.  To show that it is
also $1$-$1$ we note that a map $\cat T \to\Delta\rel X \in \RelnCat$
can be considered as a pair $(f,x)$ of maps
\begin{displaymath}
  \cat T\xrightarrow{\enspace f\enspace}\Delta\rel[-] \in \RelnCat
  \qquad\text{and}\qquad
  fT\xrightarrow{\enspace x\enspace} X \in \simp^{n}\cat S \Period
\end{displaymath}
We then have to show that if
\begin{align*}
  \cat T\xrightarrow{\enspace(f,z)\enspace}
  \Delta\rel\Delta[p_{n}, \ldots, p_{1},q]
  \qquad&\text{and}\qquad
  \Delta[p_{n}, \ldots, p_{1},q] \xrightarrow{\enspace y\enspace}
  X \qquad\text{and}\\
  \cat T\xrightarrow{\enspace (f',z')\enspace}
  \Delta\rel\Delta[p_{n}', \ldots, p_{1}',q']
  \qquad&\text{and}\qquad
  \Delta[p_{n}', \ldots, p_{1}',q'] \xrightarrow{\enspace y'\enspace} X
\end{align*}
are such that
\begin{displaymath}
  (f, yz) = (f', y'z')\colon \cat T \longrightarrow \Delta\rel X \Comma
\end{displaymath}
then these two pair represent the same element of $\colim_{\Delta
  X}\map(\cat T, \Delta\rel F)$.  This follows however from the
observation that in that case $f = f'$ and that the following diagram
commutes
\begin{displaymath}
  \xymatrix{
    {\Delta[p_{n}, \ldots, p_{1},q]} \ar[dr]_{y}
    & {fT=f'T} \ar[l]_-{z} \ar[r]^-{z'} \ar[d]
    & {\Delta[p_{n}', \ldots, p_{q}', q']} \ar[dl]^{y'}\\
    & {X}
  }
\end{displaymath}

%--------------------------------------------------------------------
\subsection{A proof of \ref{sec:PrpNTrReEq}}
\label{sec:Prf5p8}

To prove proposition~\ref{sec:PrpNTrReEq} takes more work.

The main part of the proof consists of proving (in
\ref{sec:DivCt}--\ref{sec:PrpReeMlt} below) an identical statement for
a functor
\begin{displaymath}
  K_{\delta}\colon \simp^{n}\cat S \longrightarrow \RelnCat
\end{displaymath}
and then noting (in \ref{sec:PrpCanIso}) that $K_{\delta}$ is
essentially just an alternate way of describing the functor
$\Delta\rel$.

To do this we start with considering
%--------------------------------------------------------------------
\subsection{The division of an $n$-relative category}
\label{sec:DivCt}

Given an object $\cat C \in \RelnCat$, its \textbf{division}
$\delta\cat C \in \RelnCat$ is defined as follows:
\begin{enumerate}
\item \label{Divi} $a\delta\cat C$ is the category which has as
  objects the functors $\cat p \to a\cat C$ ($p\ge 0$) and as maps
  \begin{displaymath}
    (x_{1}\colon \cat p_{1} \to a\cat C)\longrightarrow
    (x_{2}\colon \cat p_{2} \to a\cat C)
  \end{displaymath}
  the commutative diagrams of the form
  \begin{displaymath}
    \xymatrix{
      {\cat p_{1}} \ar[rr]^{f} \ar[dr]_{x_{1}}
      && {\cat p_{2}} \ar[dl]^{x_{2}}\\
      & {a\cat C}
    }
  \end{displaymath}
\end{enumerate}
and
\begin{resumeenumerate}{2}
  \item \label{Divii} $v_{i}\delta\cat C$ ($1 \le i \le n$) and
    $w\delta\cat C$ consists of those maps as in (i) for which the
    induced map
    \begin{displaymath}
      x_{1}p_{1} = x_{2}fp_{1} \longrightarrow x_{2}p_{2}
    \end{displaymath}
    is in $v_{i}\cat C$ or $w\cat C$ respectively.
\end{resumeenumerate}

Clearly $\delta\cat C$ \emph{is natural in} $\cat C$.

Moreover
\begin{resumeenumerate}{3}
  \item \label{Diviii} $\delta\cat C$ comes with a natural
    (\textbf{terminal}) \textbf{projection map}
    \begin{displaymath}
      \pi_{t}\colon \delta\cat C \longrightarrow \cat C \in \RelnCat
    \end{displaymath}
    which sends each object $x\colon \cat p \to \cat C \in \delta\cat
    C$ to the object $xp \in \cat C$ and which clearly has the
    following property:
\end{resumeenumerate}

%--------------------------------------------------------------------
\subsection{Proposition}
\label{sec:PrpfDltaiff}
\begin{em}
  A map $f \in \delta\cat C$ is in $v_{i}\delta\cat C$ ($1 \le i \le
  n$) or $w\delta\cat C$ iff $\pi_{t}f$ is in $v_{i}\cat C$ or $w\cat
  C$ respectively.
\end{em}

\intro
Using these divisions we then define
%--------------------------------------------------------------------
\subsection{A functor $K_\delta\colon \simp^n\cat S \to \RelnCat$ and
  a natural transformation $\pi_{t}\colon K_{\delta} \to K$}
\label{sec:FncNT}

We denote by
\begin{displaymath}
  K_{\delta}\colon \simp^{n}\cat S \longrightarrow \RelnCat
\end{displaymath}
the colimit preserving functor which sends each standard
$\Delta[p_{n}, \ldots, p_{1},q]$-simplex $(p_{n}, \ldots, p_{1}, q \ge
0)$ to the object
\begin{displaymath}
  \delta\cat p_{n}^{v_{n}}\times\cdots\times \delta\cat p_{1}^{v_{1}}
  \times \delta\cat q^{w} \in \RelnCat
\end{displaymath}
and with a slight abuse of notation we denote by
\begin{displaymath}
  \pi_{t}\colon K_{\delta}\longrightarrow K
\end{displaymath}
the natural transformation which is induced by the natural maps
\eqref{Diviii}
\begin{equation}
  \label{eq:starred}
  \tag{$*$}
  \delta\cat p_{n}^{v_{n}}\times \cdots\times \delta\cat
  p_{1}^{v_{1}}\times \delta\cat q^{w}
  \xrightarrow{\enspace \pi_{t}\times\cdots\times\pi_{t}\times\pi_{t}
    \enspace}
  \cat p_{n}^{v_{n}}\times\cdots\times\cat p_{1}^{v_{1}}\times\cat
  q^{w}
\end{equation}

\intro
Now we can formulate the desired \eqref{sec:Prf5p8} variation on
\ref{sec:PrpNTrReEq}:
%--------------------------------------------------------------------
\subsection{Proposition}
\label{sec:PrpReeMlt}
\begin{em}
  For every standard multisimplex $\Delta[p_{n}, \ldots, p_{1},q] \in
  \simp^{n}\cat S$, the map \eqref{sec:FncNT}
  \begin{displaymath}
    N\pi_{t}\colon NK_{\delta}\Delta[p_{n}, \ldots, p_{1},q]
    \longrightarrow
    NK\Delta[p_{n}, \ldots, p_{1},q] \in \simp^{n}\cat S
  \end{displaymath}
  is a Reedy equivalence.
\end{em}

\emsection{Proof}
As
\begin{align*}
  K_{\delta}\Delta[p_{n}, \ldots, p_{1},q] &=
  \delta\cat p_{n}^{v_{n}}\times \cdots \times \delta\cat
  p_{1}^{v_{1}} \times\delta\cat q^{w}\\
  \intertext{and}
  K\Delta[p_{n}, \ldots, p_{1},q] &=
  \cat p_{n}^{v_{n}}\times \cdots \times \cat p_{1}^{v_{1}} \times\cat
  q^{w}
\end{align*}
we have to prove that application of the functor $N$ to the map in
\ref{sec:FncNT}\eqref{eq:starred} above yields a Reedy equivalence.
But $N$ is a right adjoint and hence preserves products and it
therefore suffices to show that each of the maps
\begin{displaymath}
  N\pi_{t}\colon N\delta\cat p_{i}^{v_{i}} \longrightarrow
  N\cat p_{i}^{v_{i}}
  \qquad\text{and}\qquad
  N\pi_{t}\colon N\delta\cat a^{w} \longrightarrow N\cat q^{w}
\end{displaymath}
is a Reedy equivalence.

To do this let
\begin{displaymath}
  \tau\colon \cat p^{v_{i}} \longrightarrow \delta\cat p^{v_{i}}
  \qquad\text{and}\qquad
  \tau\colon \cat p^{w} \longrightarrow \delta\cat p^{w}
\end{displaymath}
be the maps which send an object $b \in \cat p$ to the object
\begin{displaymath}
  \cat b = (0 \to \cdots \to b) \xrightarrow{\text{ incl. }}
  (0 \to \cdots \to p) \in \delta\cat p^{v_{i}} \text{ or }\delta\cat
  p^{w} \Period
\end{displaymath}
Then $\pi_{t}\tau = 1$ and there are obvious maps
\begin{displaymath}
  h\colon \delta\cat p^{v_{i}}\times\cat 1^{w} \longrightarrow
  \delta\cat p_{v_{i}} 
  \qquad\text{and}\qquad
  h\colon \delta\cat p^{w}\times \cat 1^{w}\longrightarrow \delta\cat
  p^{w}
\end{displaymath}
such that $h0=1$ and $h1 = \tau\pi_{t}$.

The desired result then follows readily from the observation that if
\begin{enumerate}
\item \label{StrctHmi} two maps $f,g\colon \cat C \to \cat D \in
  \RelnCat$ are \emph{strictly homotopic} in the sense that there
  exists a map $h\colon \cat C\times\cat 1^{w} \to \cat D \in
  \RelnCat$ connecting them,
\end{enumerate}
then
\begin{resumeenumerate}{2}
  \item \label{StrctHmii} the maps $Nf,Ng\colon N\cat C \to N\cat D
    \in \simp^{n}\cat S$ are \emph{strictly homotopic} in the sense
    that there exists a map $k\colon N\cat C\times\Delta[0, \ldots, 0,
    1] \to N\cat D \in \simp^{n}\cat S$ connecting them,
\end{resumeenumerate}
where
\begin{resumeenumerate}{3}
  \item \label{StrctHmiii} $k$ is the composition
    \begin{multline*}
      N\cat C\times\Delta[0,\ldots, 1] \xrightarrow{\enspace \eta\enspace }
      N\cat C\times NK\Delta[0, \ldots, 0, 1]\\
      \xrightarrow{\enspace \text{Id}\enspace } N\cat C\times N\cat 1^{w}
      \xrightarrow{\enspace \iso\enspace } N(\cat C\times\cat 1^{w})
      \xrightarrow{\enspace h\enspace } N\cat D
    \end{multline*}
\end{resumeenumerate}

\intro
As mentioned in \ref{sec:Prf5p8} above,
proposition~\ref{sec:PrpNTrReEq} now is an immediate consequence of
\ref{sec:PrpReeMlt} above and
%--------------------------------------------------------------------
\subsection{Proposition}
\label{sec:PrpCanIso}
\begin{em}
  There exists a commutative diagram
  \begin{displaymath}
    \xymatrix{
      {\Delta\rel} \ar[rr] \ar[dr]
      && {K_{\delta}} \ar[dl]^{\pi_{t}}\\
      & {K}
    }
  \end{displaymath}
  of functors $\simp^{n}\cat S \to \RelnCat$ and natural
  transformations between them in which
  \begin{enumerate}
  \item the right hand map is as in \ref{sec:FncNT} and
  \item the top map is an isomorphism which, for every standard
    multisimplex $\Delta[p_{n}, \ldots, p_{1},q] \in \simp^{n}\cat S$ sends
    \begin{displaymath}
      \Delta\rel\Delta[p_{n},\ldots, p_{1},q]
      \quad\text{to}\quad
      K_{\delta}\Delta[p_{n}, \ldots, p_{1}, q] \Period
    \end{displaymath}
  \end{enumerate}
\end{em}

\emsection{Proof}
As both functors $\Delta\rel$ and $K_{\delta}$ are colimit preserving
(\ref{sec:ObvIso} and \ref{sec:FncNT}) this follows immediately from
the observation that, for every sequence of integers $p_{n}, \ldots,
p_{1}, q \ge 0$
\begin{displaymath}
  \Delta\rel\Delta[p_{n}\ldots, p_{1},q]
  \qquad\text{and}\qquad
  \delta\cat p_{n}^{v_{n}}\times\cdots \times\delta\cat
  p_{1}^{v_{1}}\times\delta\cat q^{w}
\end{displaymath}
are canonically isomorphic.

%--------------------------------------------------------------------
%--------------------------------------------------------------------
\section{Appendix}
\label{sec:App}

In this appendix we note that the categories $\RelnCat$ ($n \ge 1$)
have two additional properties which one would expect a homotopy
theory of homotopy theories to have:
\begin{lettered}
  \item \label{AppxA} There there exists a functor $\RelnCat\to
    \RelnpoCat$ which has a \emph{left inverse right adjoint}.
  \item \label{AppxB} That every object of $\RelnpoCat$ gives rise to
    a \emph{category enriched over $\RelnCat$} which suggests the
    possibility that ``a map in $\RelnpoCat$ is a Rezk equivalence
    \eqref{sec:ReeRzEq} iff the induced map between these enriched
    categories is a kind of $DK$-equivalence''.
\end{lettered}

\intro
To deal with \ref{AppxA} we note that a straightforward calculation
yields:
%--------------------------------------------------------------------
\subsection{Proposition}
\label{sec:PrpRtAdjInv}
\begin{em}
  For every integer $n \ge 1$ the functor
  \begin{displaymath}
    \RelnCat \longrightarrow \RelnpoCat
  \end{displaymath}
  which sends
  \begin{displaymath}
    (a\cat C, v_{1}\cat C, \ldots, v_{n}\cat C, w\cat C)
    \qquad\text{to}\qquad
    (a\cat C, v_{1}\cat C, \ldots, v_{n}\cat C, w\cat C, w\cat C)
  \end{displaymath}
  has a right adjoint left inverse which sends
  \begin{displaymath}
    (a\cat D, v_{1}\cat D, \ldots, v_{n+1}\cat D, w\cat D)
    \qquad\text{to}\qquad
    (\bar a\cat D, v_{1}\cat D, \ldots, v_{n}\cat D, w\cat D)
  \end{displaymath}
  where $\bar a\cat D \subset a\cat D$ denotes the subcategory which
  consists of the finite compositions of maps in the $v_{i}\cat D$ ($1
  \le i \le n$).
\end{em}

\intro
We deal with \ref{AppxB} by means of an $n$-relative version of the
Grothendieck enrichment of \cite{DHKS}*{3.4 and 3.5}.

To do this we start with recalling
%--------------------------------------------------------------------
\subsection{Types of zigzags}
\label{sec:ZgZg}

The \textbf{type} of a zigzag of maps in a category $\cat C$ from an
object $X$ to an object $Y$
\begin{displaymath}
  \xymatrix{
    {X} \ar@{-}[r]^{f_{1}}
    & {\cdot} \ar@{}[r]|{\spacedcdots}
    & {\cdot} \ar@{-}[r]^{f_{m}}
    & {Y}
    & {(m \ge 0)}
  }
\end{displaymath}
will be the pair $T = (T_{+}, T_{-})$ of complementary subsets of the
set of integers $\{1, \ldots, m\}$ such that $i \in T_{+}$ whenever
$f_{i}$ is a forward map and $i \in T_{-}$ otherwise.

These types can be considered as the objects of a \textbf{category of
  types} $\cat T$ which has, for every two types $(T_{+}, T_{-})$ and
$(T_{+}', T_{-}')$ of length $m$ and $m'$ respectively, as maps
$t\colon (T_{+}, T_{-}) \to (T_{+}', T_{-}')$ the weakly monotonic
maps $t\colon \{1, \ldots, m\} \to \{1, \ldots, m'\}$ such that
\begin{displaymath}
  tT_{+} \subset T_{+}'
  \qquad\text{and}\qquad
  tT_{-} \subset T_{-}' \Period
\end{displaymath}

\intro
With these types one then associates
%--------------------------------------------------------------------
\subsection{$n$-Relative arrow categories}
\label{sec:nRlArCt}

Given an object $\cat C \in \RelnpoCat$ let, as in
\ref{sec:PrpRtAdjInv}, $\bar a\cat C \subset a\cat C$ denote the
subcategory which consists of the finite compositions of maps of the
$v_{i}\cat C$ ($1 \le i \le n$).

For every pair of objects $X,Y \in \cat C$ and type $T$
\eqref{sec:ZgZg} we then denote by $\cat C^{T}(X,Y) \in \RelnCat$ the
\textbf{$n$-relative arrow category} which has
\begin{enumerate}
\item \label{NRlAri} as objects the zigzags of type $T$ in $\cat C$
  between $X$ and $Y$ in which the backward maps are in $\bar a\cat
  C$,
\item \label{NRlArii} as maps in $v_{i}\cat C^{T}(X,Y)$ ($1 \le i \le
  n$) and $w\cat C^{T}(X,Y)$ between two such zigzags the commutative
  diagrams of the form
  \begin{displaymath}
    \xymatrix{
      {X} \ar@{-}[r] \ar[d]_{1}
      & {\cdot} \ar@{}[r]|{\spacedcdots} \ar[d]
      & {\cdot} \ar@{-}[r] \ar[d]
      & {Y} \ar[d]^{1}\\
      {X} \ar@{-}[r]
      & {\cdot} \ar@{}[r]|{\spacedcdots}
      & {\cdot} \ar@{-}[r]
      & {Y}
    }
  \end{displaymath}
  in which the vertical maps are in $v_{i}\cat C$ and $w\cat C$
  respectively, and
\item \label{NRlAriii} as maps in $a\cat C^{T}(X,Y)$ the finite
  compositions of maps of the $v_{i}\cat C^{T}(X,Y)$ ($1 \le i \le
  n$).
\end{enumerate}

\intro
These arrow categories in turn give rise to
%--------------------------------------------------------------------
\subsection{$\cat T$-diagrams of arrow categories}
\label{sec:TDiagArCt}

Given an object $\cat C \in \RelnCat$ and objects $X,Y \in \cat C$,
one can form a \textbf{$\cat T$-diagram of arrow categories}
\begin{displaymath}
  \cat C^{(\cat T)}(X,Y)\colon \cat T \longrightarrow \RelnCat
\end{displaymath}
which assigns to every object $T \in \cat T$ the arrow category
\begin{displaymath}
  \cat C^{T}(X,Y) \in \RelnCat
\end{displaymath}
and to every map $t\colon T \to T' \in \cat T$ the map
\begin{displaymath}
  t_{*}\colon \cat C^{T}(X,Y)\longrightarrow \cat C^{T'}(X,Y)
  \in \RelnCat
\end{displaymath}
which sends a zigzag of type $T$
\begin{displaymath}
  \xymatrix{
    {X} \ar@{-}[r]^{f_{1}}
    & {\cdot} \ar@{}[r]|{\spacedcdots}
    & {\cdot} \ar@{-}[r]^{f_{m}}
    & {Y}
  }
\end{displaymath}
to the zigzag of type $T'$
\begin{displaymath}
  \xymatrix{
    {X} \ar@{-}[r]^{f_{1}'}
    & {\cdot} \ar@{}[r]|{\spacedcdots}
    & {\cdot} \ar@{-}[r]^{f_{m}'}
    & {Y}
  }
\end{displaymath}
in which each $f_{j}'$ ($1 \le j \le m'$) is the composition of the
$f_{i}$ with $ti=j$ or, in no such $i$ exists, the appropriate
identity map.

\intro
Now we can form
%--------------------------------------------------------------------
\subsection{The Grothendieck construction on $\cat C^{(\cat T)}(X,Y)$}
\label{sec:GrthCT}

Given an object $\cat C \in \RelnpoCat$ and objects $X,Y \in \cat C$
the \textbf{Grothendieck construction} on $\cat C^{(\cat T)}(X,Y)$ is
the object
\begin{displaymath}
  \Gr\cat C^{(\cat T)}(X,Y) \in \RelnCat
\end{displaymath}
which has
\begin{enumerate}
\item \label{Gri} as objects the zigzags in $\cat C$ between $X$ and
  $Y$ in which the backward maps are in $\bar a\cat C$, i.e.\ pairs
  $(T,Z)$ consisting of objects
  \begin{displaymath}
    T\in\cat T
    \qquad\text{and}\qquad
    Z \in \cat C^{T}(X,Y)
  \end{displaymath}
\end{enumerate}
and
\begin{resumeenumerate}{2}
  \item \label{Grii} for every two such objects $(T,Z)$ and $(T',
    Z')$, as maps $(T,Z) \to (T', Z')$ the pairs $(t,z)$ consisting of
    maps
    \begin{displaymath}
      t\colon T \longrightarrow T' \in \cat T
      \qquad\text{and}\qquad
      z\colon t_{*}Z \longrightarrow Z' \in \cat C^{T'}(X,Y)
    \end{displaymath}
\end{resumeenumerate}
and in which
\begin{resumeenumerate}{3}
  \item \label{Griii} for every two composable maps $(t,z)$ and
    $(t',z')$ their composition is defined by the formula
    \begin{displaymath}
      (t',z')(t,z) = \bigl(t't, z'(t_{*}z)\bigr)
    \end{displaymath}
\end{resumeenumerate}

\intro
Together these Grothendieck constructions give rise to
%--------------------------------------------------------------------
\subsection{A Grothendieck enrichment}
\label{sec:GrEnrch}

Given an object $\cat C \in \RelnpoCat$ we now define its
\textbf{Grothendieck enrichment} as the \emph{category $\Gr\cat
  C^{(\cat T)}$ enriched over $\RelnCat$} which
\begin{enumerate}
\item \label{GrEni} has the same objects as $\cat C$,
\item \label{GrEnii} has for every two objects $X,Y \in \cat C$, as
  it's hom-object the $n$-relative category $\cat C^{(\cat T)}(X,Y)$, and
\item \label{GrEniii} has, for every three objects $X$, $Y$ and $Z \in
  \cat C$ as composition
  \begin{displaymath}
    \Gr\cat C^{(\cat T)}(X,Y) \times \Gr\cat C^{(\cat T)}(X,Y)
    \longrightarrow \Gr\cat C^{(\cat T)}(X,Z)
  \end{displaymath}
  the function induced by the compositions of the zigzags involved.
\end{enumerate}

%--------------------------------------------------------------------
%--------------------------------------------------------------------

\end{document}